\input amstex
\documentstyle{amsppt}
\pagewidth{5.4in}
\pageheight{7.6in}
\magnification=1200
\TagsOnRight
\NoRunningHeads
\topmatter
\title
\bf A simple proof on the non-existence of\\ 
shrinking breathers for the Ricci flow 
\endtitle
\author
Shu-Yu Hsu
\endauthor
\affil
Department of Mathematics\\
National Chung Cheng University\\
168 University Road, Min-Hsiung\\
Chia-Yi 621, Taiwan, R.O.C.\\
e-mail:syhsu\@math.ccu.edu.tw
\endaffil
\date
Jan 16, 2006
\enddate
\address
e-mail address:syhsu\@math.ccu.edu.tw
\endaddress
\abstract
Suppose $M$ is a compact n-dimensional manifold, $n\ge 2$, with a metric 
$g_{ij}(x,t)$ 
that evolves by the Ricci flow $\partial_tg_{ij}=-2R_{ij}$ in $M\times (0,T)$. 
We will give a simple proof of a recent result of Perelman on the 
non-existence of shrinking breather without using the logarithmic Sobolev 
inequality.  
\endabstract
\endtopmatter
\NoBlackBoxes
\define \pd#1#2{\frac{\partial #1}{\partial #2}}
\define \1{\partial}
\define \2{\overline}
\define \3{\varepsilon}
\define \4{\widetilde}
\define \5{\underline}
\document

It is known that Ricci flow is a very powerful tool in understanding the 
geometry and structure of manifolds. In 1982 R.~Hamilton \cite{H1} first 
began the study of Ricci flow on a manifold.
Suppose $M$ is a compact 3-dimensional manifold with a metric $g_{ij}(x)$ 
having a strictly positive Ricci curvature. R.~Hamilton proved that if  
the metric $g_{ij}$ evolve by the Ricci flow
$$
\frac{\1 }{\1 t}g_{ij}=-2R_{ij}\tag 0.1
$$
with $g_{ij}(x,0)=g_{ij}(x)$, then the evolving metric will converge 
modulo scaling to a metric of constant positive curvature. A similar result
for compact 4-dimensional manifold with positive curvature operator was
proved by R.~Hamilton in the paper \cite{H2}. By using a modification of the
proof of Li-Yau Harnack inequality \cite{LY} for the heat equations on 
manifolds R.~Hamilton \cite{H4} proved the Harnack inequality for the Ricci 
flow. Singularities of solutions of the Ricci flow was studied by R.~Hamilton
\cite{H5} and G.~Perelman \cite{P1}, \cite{P2}.

Ricci flow on non-compact manifolds was studied by W.X.~Shi \cite{S1},
\cite{S2}, R.~Hamilton \cite{H3}, and L.F.~Wu \cite{W1}, \cite{W2}. 
Existence and asymptotic behaviour of solutions of the Ricci flow equation on 
non-compact $\Bbb{R}^2$ was studied by S.Y.~Hsu in the papers \cite{Hs1--4}. 
We refer the reader to the paper \cite{H5} by R.~Hamilton and 
the book \cite{CK} B.~Chow and D.~Knopf for various recent results on the 
Ricci flow. One can also read the recent lecture notes by B.~Chow \cite{C}
on Ricci flow.

A metric $g_{ij}(t)$ evolving by the Ricci flow in $M\times (0,T)$ is called 
a steady (shrinking, expanding respectively) breather if there exist 
$0<t_1<t_2<T$ and 
$\alpha =1$ ($0<\alpha<1$, $\alpha>1$ respectively) and a diffeomorphism 
$\phi :M\to M$ such that $g_{ij}(t_2)=\phi^{\ast}(\alpha g_{ij}(t_1))$.
As observed by G.~Perelman \cite{P1} if one considers Ricci flow as a dynamical
system on the space of Riemannian metrics modulo diffeomorphism and scaling, 
then breathers correspond to periodic orbits for the Ricci flow. So it is 
interesting to know whether breather exists in a Ricci flow.

In the paper \cite{P1} G.~Perelman found two functionals for the Ricci flow 
which are monotone increasing with respect to time.  G.~Perelman then used 
these and logarithmic Sobolev inequality to proved that there is no expanding 
or shrinking breathers for the Ricci flow. However his proof of non-existence 
of shrinking breathers has some gaps and requires the existence of solution 
of some auxillary parabolic equation on a manifold with initial value a delta 
mass which is highly non-trivial. In this paper we will modify Perelman's 
argument and give a simple proof of the non-existence of shrinking breathers
without using the logarithmic Sobolev inequality. 

The plan of the paper is as follows. In section 1 we will prove some technical
lemmas. In section 2 we will fix the gaps in the proof of the monotonicity
property of the $W$ functional in Perelman's paper \cite{P1}. We will also
prove the non-existence of shrinking breathers.

We will assume that $M$ is a compact n-dimensional manifold, $n\ge 2$, with 
a metric 
$g(t)=(g_{ij}(\cdot,t))$ that evolves by the Ricci flow (0.1) in $M\times 
(0,T)$ for the rest of the paper. 

$$
\text{Section 1}
$$

In this section we will establish some technical lemmas. We first recall a 
standard result (cf. Theorem 1.6.2 of \cite{J}).

\proclaim{\bf Lemma 1.1}
Let $0<t_1<T$ and $f\in C^{\infty}(M\times (0,t_1))$. For any 
$t\in (0,t_1)$ there exist a smooth function $\psi_s^t(p)
=\psi^t(p,s):M\times (0,t_1)\to M$
satisfying
$$\left\{\aligned
\frac{\1}{\1 s}\psi^t(p,s)=&-\nabla f(\psi^t(p,s),t)\quad\forall 
s\in (0,t_1),p\in M\\
\psi^t(p,0)=&p\qquad\qquad\qquad\qquad\forall p\in M.\endaligned\right.\tag 1.1
$$
\endproclaim

By an argument similar to the proof of Theorem 1.6.2 of \cite{J} we have the 
following lemma.

\proclaim{\bf Lemma 1.2}
Let $0<t_1<T$ and $f\in C^{\infty}(M\times (0,t_1))$. For any 
$t_0\in (0,t_1)$ there exist a smooth function $\phi_{t_0}(p,t)=\phi_{t_0,t}
(p)$ such that $\phi_{t_0}:M\times (0,t_1)\to M$ and  satisfies
$$\left\{\aligned
\frac{\1}{\1 t}\phi_{t_0}(p,t)=&-\nabla f(\phi_{t_0}(p,t),t)\quad\forall 
t\in (0,t_1),p\in M\\
\phi_{t_0}(p,t_0)=&p\qquad\qquad\qquad\qquad\forall p\in M.
\endaligned\right.\tag 1.2
$$
If $t_0'\in (0,t_1)$, then the map $\phi_{t_0,t_0'}:M\to M$
is a diffeomorphism with inverse $\phi_{t_0',t_0}$.
\endproclaim

\proclaim{\bf Lemma 1.3}
Let $0<t_1<T$, $t,t_0\in (0,t_1)$, $\delta_0=\min (t,t_1-t)$, and let $f$, 
$\psi_s^{t}$, $\phi_{t_0,t}$ be as in Lemma 1.1 and Lemma 1.2. Let $p\in M$ 
and $x=(x_1,\dots,x_n):U\subset M\to\Bbb{R}^n$ be a local co-ordinate chart 
around $p_0=\phi_{t_0,t}(p)$ for some open neighbourhood $U$ of $p_0$ 
such that $x(U)=B_{R_0}$ for some $R_0>0$ and $x(\phi_{t,t+h}(p_0)), 
x(\psi_h^t(p_0))\in B_{R_0}$ for any $|h|\le\delta_1$ for 
some constant $0<\delta_1\le\delta_0$. Let $e(h)
=x(\phi_{t,t+h}(p_0))-x(\psi_h^{t}(p_0))$ for any $|h|\le\delta_1$. Then 
there exists a constant $C>0$ such that
$$
|e(h)|+\max_{1\le k\le n}\left |\frac{d e_k}{dh}\right |
+\max_{1\le j,k\le n}\biggl\{\left |\biggl (\frac{\1}{\1 x_j}
\biggr )_{p_0}e_k\right |
+\left |\frac{d }{d h}\biggl (\frac{\1}{\1 x_j}
\biggr )_{p_0}e_k\right |\biggr\}
\le C|h|\qquad\forall |h|\le\delta_1\tag 1.3
$$
where $e(h)=(e_k(h))_{k=1}^n$ in this local co-ordinate system
and $|e(h)|=(\sum_{k=1}^ne_k(h)^2)^{1/2}$. 
\endproclaim
\demo{Proof}
Without loss of generality we will abuse the notation and write 
$\phi_{t,t+h}(p_0)$, $\psi_h^{t}(p_0)$, instead of $x(\phi_{t,t+h}(p_0))$, 
$x(\psi_h^{t}(p_0))$, etc. and we will write $\1/\1 x_j$ for 
$(\1/\1 x_j)_{p_0}$. Let $\phi_{t,t+h}(p_0)
=(\phi_{t,t+h}^k(p_0))_{k=1}^n$ and $\psi_h^{t}(p_0)
=((\psi_h^{t})^k(p_0))_{k=1}^n$ in the local co-ordinate system $(x,U)$ and
let
$$
q(s)=s\phi_{t,t+h}(p_0)+(1-s)\psi_h^{t}(p_0)\quad\forall 0\le s\le 1.\tag 1.4
$$
By (1.1) and (1.2),
$$\align
&\biggl |\frac{d e_k}{dh}\biggr |\\
=&\biggl |\frac{d}{d h}\phi_{t,t+h}^k(p_0)
-\frac{d}{dh}(\psi_h^{t})^k(p_0)\biggr |\\
=&\biggl |-g^{kj}(\phi_{t,t+h}(p_0),t+h)\frac{\1}{\1 x_j}f
(\phi_{t,t+h}(p_0),t+h)\\
&\qquad +g^{kj}(\psi_h^{t}(p_0),t)
\frac{\1}{\1 x_j}f(\psi_h^{t}(p_0),t)\biggr |\\
\le&\biggl |g^{kj}(\phi_{t,t+h}(p_0),t)\frac{\1}{\1 x_j}f(\phi_{t,t+h}(p_0),t)
-g^{kj}(\phi_{t,t+h}(p_0),t+h)\frac{\1}{\1 x_j}f(\phi_{t,t+h}(p_0),t+h)
\biggr |\\
&\quad +\biggl |g^{kj}(\psi_h^{t}(p_0),t)
\frac{\1}{\1 x_j}f(\psi_h^{t}(p_0),t)
-g^{kj}(\phi_{t,t+h}(p_0),t)\frac{\1}{\1 x_j}f(\phi_{t,t+h}(p_0),t)\biggr |\\
\le&C|h|+\biggl |\int_0^1\frac{d}{d s}\biggl (g^{kj}(q(s),t)\frac{\1}{\1 x_j}
f(q(s),t)\biggr )\,ds\biggr |\\
\le&C|h|+\biggl |\int_0^1e_i\biggl (\frac{\1}{\1 x_i}g^{kj}(q(s),t)
\frac{\1}{\1 x_j}f(q(s),t)
+g^{kj}(q(s),t)\frac{\1^2}{\1 x_i\1 x_j}f(q(s),t)\biggr )
\,ds\biggr |\\
\le& C(|h|+|e|)\qquad\qquad\qquad\qquad\forall |h|\le\delta_1,k=1,2,\dots,n.
\tag 1.5\endalign
$$
Hence
$$\align
&\biggl |\frac{d |e|^2}{dh}\biggr |\le C(|e|^2+h^2)\qquad\qquad\,\,
\forall  |h|\le\delta_1\\
\Rightarrow\quad&\biggl |\frac{d}{dh}( e^{-Ch}|e(h)|^2)\biggr |\le Ch^2 e^{-Ch}
\quad\forall  |h|\le\delta_1\\
\Rightarrow\quad&|e(h)|^2\le C'h^2\qquad\qquad\qquad\qquad
\forall  |h|\le\delta_1.\tag 1.6
\endalign
$$
Similarly
$$
\left |\frac{d}{d h}\biggl (\frac{\1 e_k}{\1 x_j}\biggr )\right |
\le C(|h|+|e|)\le C|h|
\quad\Rightarrow\quad\left |\frac{\1 e_k}{\1 x_j}\right |\le C|h|\quad\forall
|h|\le\delta_1,j,k=1,2,\dots,n.\tag 1.7
$$
By (1.5), (1.6), and (1.7) we get (1.3) and the lemma follows.
\enddemo

\proclaim{\bf Lemma 1.4}
Let $0<t_1<T$, $t_0\in (0,t_1)$, and let $f$, $\phi_{t_0,t}$, be as in 
Lemma 1.2. Let 
$$
\2{g}(t)=\phi_{t_0,t}^{\ast}(g(t))\quad\forall 0<t<t_1.\tag 1.8
$$
Then
$$
\frac{\1}{\1 t}\2{g}(t)=\phi_{t_0,t}^{\ast}\biggl (\frac{\1}{\1 t}g(t)+
L_{V(t)}(g(t))\biggr )\quad\forall 0<t<t_1\tag 1.9
$$
where $V(t)=-\nabla f(\cdot,t)$.
\endproclaim
\demo{Proof}
Let $p\in M$, $t\in (0,t_1)$, and let $\psi_s^t$ be as in Lemma 1.1. Let 
$(x,U)$, $\delta_1>0$, and $e(h)=(e_k(h))_{k=1}^n$ be as in 
Lemma 1.3. Let $X$ and $Y$ be two vector fields on $M$. Then there exist
a constant $\delta_2\in (0,\delta_1)$ and an open neighbourhood 
$V\subset U$ of the curve $t'\to\phi_{t_0,t'}(p)$, $t-\delta_2\le t'\le t
+\delta_2$, such that $x(V)$ is convex in $\Bbb{R}^n$ and the vector fields 
$d\phi_{t_0,t'}(X(p))$ and $d\phi_{t_0,t'}(Y(p))$ along the curve $t'\to
\phi_{t_0,t'}(p)$, $t-\delta_2\le t'\le t+\delta_2$, can be extended to two 
local vector fields $\4{X}$ and $\4{Y}$ on $V$. That is
$$\left\{\aligned
&\4{X}(\phi_{t_0,t'}(p))=d\phi_{t_0,t'}(X(p))\quad\forall t'\in (t-\delta_2,
t+\delta_2)\\
&\4{Y}(\phi_{t_0,t'}(p))=d\phi_{t_0,t'}(Y(p))\quad\forall t'\in (t-\delta_2,
t+\delta_2).\endaligned\right.
$$
Let $p_0=\phi_{t_0,t}(p)$ and
$$\align
E(h)=&g(\phi_{t_0,t+h}(p),t)(d\phi_{t_0,t+h}(X(p)),d\phi_{t_0,t+h}(Y(p)))\\
&\qquad -g(\psi_h^t(p_0),t)(d\psi_h^t(\4{X}(p_0)),d\psi_h^t(\4{Y}(p_0))).
\endalign
$$ 
Let $\phi_{t,t+h}=(\phi_{t,t+h}^k)_{k=1}^n$ and $\psi_h^{t}
=((\psi_h^{t})^k)_{k=1}^n$  in the local co-ordinate system $(x,U)$. We write
$$\left\{\aligned
&\4{X}(q)=a^i(q)\frac{\1}{\1 x_i}\\
&\4{Y}(q)=b^i(q)\frac{\1}{\1 x_i}
\endaligned\right.
$$
and let $q(s)=(q(s)^k)_{k=1}^n$ be given by (1.4). Since $\phi_{t_0,t+h}
=\phi_{t,t+h}\circ\phi_{t_0,t}$ on $M$, 
$$\align
E(h)=&g(\phi_{t,t+h}(p_0),t)(d\phi_{t,t+h}(\4{X}(p_0)),
d\phi_{t,t+h}(\4{Y}(p_0)))\\
&\qquad -g(\psi_h^{t}(p_0),t)(d\psi_h^t(\4{X}(p_0)),d\psi_h^t(\4{Y}(p_0)))\\
=&g_{ij}(\phi_{t,t+h}(p_0),t)\frac{\1\phi_{t,t+h}^i}{\1 x_k}(p_0)
\frac{\1\phi_{t,t+h}^j}{\1 x_{k'}}(p_0)a^k(p_0)b^{k'}(p_0)\\
&\qquad-g_{ij}(\psi_h^t(p_0),t)\frac{\1 (\psi_h^t)^i}{\1 x_k}(p_0)
\frac{\1 (\psi_h^t)^j}{\1 x_{k'}}(p_0)a^k(p_0)b^{k'}(p_0)\\
=&\int_0^1\frac{d}{ds}\biggl [g_{ij}(q(s),t)\frac{\1 q(s)^i}{\1 x_k}
\frac{\1 q(s)^j}{\1 x_{k'}}a^k(p_0)b^{k'}(p_0)\biggr ]ds\\
=&e_l(h)\int_0^1\frac{\1 g_{ij}}{\1 x_l}(q(s),t)
\frac{\1 q(s)^i}{\1 x_k}\frac{\1 q(s)^j}{\1 x_{k'}}a^k(p_0)b^{k'}(p_0)\,ds\\
&\qquad
+\int_0^1g_{ij}(q(s),t)\frac{\1 e_i}{\1 x_k}\frac{\1 q(s)^j}{\1 x_{k'}}a^k
(p_0)b^{k'}(p_0)\,ds\\
&\qquad
+\int_0^1g_{ij}(q(s),t)\frac{\1 q(s)^i}{\1 x_k}\frac{\1 e_j}{\1 x_{k'}}
a^k(p_0)b^{k'}(p_0)\,ds\\
=&E_1(h)+E_2(h)+E_3(h)\qquad\qquad\qquad\quad\forall |h|\le\delta_2.\tag 1.10
\endalign
$$
Let 
$$\left\{\aligned
&G_l=\frac{\1 g_{ij}}{\1 x_l}(q(s),t)\frac{\1 q(s)^i}{\1 x_k}
\frac{\1 q(s)^j}{\1 x_{k'}}a^k(p_0)b^{k'}(p_0)\quad\forall l=1,2,\dots,n\\
&H_i=g_{ij}(q(s),t)\frac{\1 q(s)^j}{\1 x_{k'}}b^{k'}(p_0)
\qquad\qquad\qquad\quad\forall i=1,2,\dots,n.
\endaligned\right.
$$
Then
$$\align
\frac{d}{d h}E_1(h)=&e_l(h)\int_0^1\biggl (s\frac{d}{d h}\phi_{t,t+h}^m(p_0)
+(1-s)\frac{d}{d h}(\psi_h^t)^m(p_0)\biggr )
\frac{\1 G_l}{\1 x_m}\,ds\\
&\qquad +\biggl (\frac{d}{d h}e_l(h)\biggr )\int_0^1G_l\,ds\\
=&-e_l(h)\int_0^1\biggl (s\nabla_m f(\phi_{t,t+h}(p_0),t+h)
+(1-s)\nabla_m f((\psi_h^t)(p_0),t)\biggr )
\frac{\1 G_l}{\1 x_m}\,ds\\
&\qquad +\biggl (\frac{d}{d h}e_l(h)\biggr )\int_0^1G_l\,ds\endalign
$$
and
$$\align
&\frac{d}{d h}E_2(h)\\
=&\int_0^1\biggl (s\frac{d}{d h}\phi_{t,t+h}^m(p_0)
+(1-s)\frac{d}{d h}(\psi_h^t)^m(p_0)\biggr )
\frac{\1 H_i}{\1 x_m}\frac{\1 e_i}{\1 x_k}a^k(p_0)\,ds\\
&\qquad+\int_0^1H_i\frac{d }{dh}\biggl (\frac{\1 e_i}{\1 x_k}
\biggr )a^k(p_0)\,ds\\
=&-\int_0^1\biggl (s\nabla_m f(\phi_{t,t+h}(p_0),t+h)
+(1-s)\nabla_m f((\psi_h^t)(p_0),t)\biggr )
\frac{\1 H_i}{\1 x_m}\frac{\1 e_i}{\1 x_k}a^k(p_0)\,ds\\
&\qquad+\int_0^1H_i\frac{d }{dh}\biggl (\frac{\1 e_i}{\1 x_k}
\biggr )a^k(p_0)\,ds.
\endalign
$$
Hence by Lemma 1.3,
$$\align
&\biggl |\frac{d}{d h}E_1(h)\biggr |+\biggl |\frac{d}{d h}E_2(h)\biggr |\\
\le& C\biggl\{|e(h)|+\max_{1\le k\le n}\left |\frac{d e_k}{dh}\right |
+\max_{1\le j,k\le n}\biggl [\left |\frac{\1 e_k}{\1 x_j}\right |
+\left |\frac{d }{d h}\biggl (\frac{\1 e_k}{\1 x_j}\biggr )\right |\biggr ]
\biggr\}\le C|h|\quad\forall |h|\le\delta_2\\
\Rightarrow\quad&\frac{d}{d h}E_1(h)\biggr |_{h=0}
=\frac{d}{d h}E_2(h)\biggr |_{h=0}=0.\tag 1.11
\endalign
$$
Similarly,
$$
\frac{d}{d h}E_3(h)\biggr |_{h=0}=0.\tag 1.12
$$
By (1.10), (1.11), and (1.12),
$$
\frac{d}{d h}E(h)\biggr |_{h=0}=0.\tag 1.13
$$
Then
$$\align
\frac{d}{d t}\2{g}(p,t)(X,Y)=&\frac{d}{d h}\biggr |_{h=0}\phi_{t_0,t+h}^{\ast}
(g(t+h))(p)(X,Y)\\
=&\frac{d}{d h}\biggr |_{h=0}g(\phi_{t_0,t+h}(p),t+h)
(d\phi_{t_0,t+h}(X(p)),d\phi_{t_0,t+h}(Y(p)))\\
=&\frac{d}{d h}\biggr |_{h=0}g(\phi_{t_0,t+h}(p),t)
(d\phi_{t_0,t+h}(X(p)),d\phi_{t_0,t+h}(Y(p)))\\
&\qquad +\frac{\1}{\1 h}\biggr |_{h=0}g(\phi_{t_0,t}(p),t+h)
(d\phi_{t_0,t}(X(p)),d\phi_{t_0,t}(Y(p)))\\
=&\frac{d}{d h}\biggr |_{h=0}g(\psi_h^{t}\circ\phi_{t_0,t}(p),t)
(d\psi_h^{t}(\4{X}(\phi_{t_0,t}(p))),d\psi_h^{t}(\4{Y}(\phi_{t_0,t}(p))))\\
&\qquad +\phi_{t_0,t}^{\ast}\biggl (\frac{\1}{\1 t}g(t)\biggr )(p)(X,Y)
\qquad\qquad\qquad\qquad\qquad\qquad(\text{by} (1.13))\\
=&L_{V(t)}(g(t))(\phi_{t_0,t}(p))(\4{X},\4{Y})
+\phi_{t_0,t}^{\ast}\biggl (\frac{\1}{\1 t}g(t)\biggr )(p)(X,Y)\\
=&L_{V(t)}(g(t))(\phi_{t_0,t}(p))(d\phi_{t_0,t}(X),d\phi_{t_0,t}(Y))
+\phi_{t_0,t}^{\ast}\biggl (\frac{\1}{\1 t}g(t)\biggr )(p)(X,Y)\\
=&\phi_{t_0,t}^{\ast}(L_{V(t)}(g(t)))(p)(X,Y)
+\phi_{t_0,t}^{\ast}\biggl (\frac{\1}{\1 t}g(t)\biggr )(p)(X,Y)\\
=&\phi_{t_0,t}^{\ast}\biggl (\frac{\1}{\1 t}g(t)+L_{V(t)}(g(t))\biggr )(p)(X,Y)
\endalign
$$
and (1.9) follows.
\enddemo

$$
\text{Section 2}
$$

In this section we will modify Perelman's argument \cite{P1} and give a 
simple proof of a recent result of Perelman \cite{P1} on the non-existence of 
shrinking breather without using the logarithmic Sobolev inequality
and the questionable existence of solution of some parabolic equation
with initial data a delta mass.

Similar to \cite{P1} for any $\tau>0$, $f\in C^{\infty}(M)$, and 
Riemannian metric $\4{g}=(\4{g}_{ij})$ on $M$, let 
$$
\Cal{F}(\4{g},f)=\int_M(R(\4{g})+|\nabla f|^2)e^{-f}\,dV_{\4{g}},
\tag 2.1
$$
$$
W(\4{g},f,\tau)=(4\pi\tau)^{-n/2}\int_M\{\tau (R(\4{g})+|\nabla f|^2)
+f-n\}e^{-f}\,dV_{\4{g}},\tag 2.2
$$
and 
$$
\mu (\4{g},\tau)=\inf_{f\in\Cal{A}(\4{g},\tau)}W(\4{g},f,\tau)\tag 2.3
$$
where $R(\4{g})$ is the scalar curvature of $\4{g}$ and
$$
\Cal{A}(\4{g},\tau)=\biggl\{f\in C^{\infty}(M):(4\pi\tau)^{-n/2}\int_Me^{-f}\,
dV_{\4{g}}=1\biggr\}.\tag 2.4
$$
We will first prove that $\mu (\4{g},\tau)$ is well-defined.

\proclaim{\bf Lemma 2.1}
Let $\4{g}$ be a Riemannian metric on $M$ and $\tau_0>0$. Then there exist 
constants $0<\delta<1$, $C_1>0$, and $C_{\tau_0}>0$ such that
$$
W(\4{g},f,\tau)
\ge\biggl [(1-\delta)\lambda_1
-\frac{\delta}{(4\pi\tau)^{\frac{n}{2}}}
\|R(\4{g})\|_{L^{\infty}(M)}-4\delta\biggr ]\tau
-C_{\tau_0}-\frac{C_1}{\tau^{n/2}}
\log(4\pi\tau)\tag 2.5 
$$
holds for any $\tau\ge\tau_0,f\in\Cal{A}(\4{g},\tau)$, where 
$\lambda_1$ is the first eigenvalue of the operator 
$R(\4{g})-4\Delta_{\4{g}}$. Hence $\mu (\4{g},\tau)>-\infty$ is well-defined
for any $\tau>0$.
\endproclaim
\demo{Proof}
Without loss of generality we may assume that $n\ge 3$.
Let $f\in\Cal{A}(\4{g},\tau)$, $\Phi =e^{-f/2}$, and 
$$
\2{W}(\4{g},\Phi,\tau)=(4\pi\tau)^{-n/2}\int_M\{\tau (R(\4{g})\Phi^2
+4|\nabla\Phi|^2)-\Phi^2\log\Phi^2\}\,dV_{\4{g}}-n.
$$
Then
$$
W(\4{g},f,\tau)=\2{W}(\4{g},\Phi,\tau)\tag 2.6
$$
and
$$
(4\pi\tau)^{-n/2}\int_M\Phi ^2\,dV_{\4{g}}=1.\tag 2.7
$$
Let $\delta\in (0,1)$. By (2.6) and (2.7),
$$\align
&W(\4{g},f,\tau)\\
=&\tau\frac{\int_M(R(\4{g})\Phi^2	
+4|\nabla\Phi|^2)\,dV_{\4{g}}}{\int_M\Phi^2\,dV_{\4{g}}}
-\frac{\int_M\Phi^2\log\Phi^2\,dV_{\4{g}}}{\int_M\Phi ^2\,dV_{\4{g}}}-n\\
\ge&(1-\delta)\tau\cdot\inf\Sb\psi\in C^{\infty}(M)\\\psi\ne 0\endSb
\biggl (\frac{\int_M(R(\4{g})\psi^2
+4|\nabla\psi|^2)\,dV_{\4{g}}}{\int_M\psi^2\,dV_{\4{g}}}\biggr )
-\frac{\delta\tau}{(4\pi\tau)^{\frac{n}{2}}}\|R(\4{g})\|_{L^{\infty}(M)}
+I(\Phi)-n\\
\ge&(1-\delta)\tau\lambda_1
-\frac{\delta\tau}{(4\pi\tau)^{\frac{n}{2}}}\|R(\4{g})\|_{L^{\infty}(M)}
+I(\Phi)-n\tag 2.8
\endalign
$$
where $\lambda_1$ is the first eigenvalue of $R(\4{g})-4\Delta_{\4{g}}$
and
$$
I(\Phi)=\frac{4\delta\tau\int_M|\nabla\Phi|^2\,dV_{\4{g}}
-\int_M\Phi^2\log\Phi^2\,dV_{\4{g}}}{\int_M\Phi ^2\,dV_{\4{g}}}.
$$
We will now use a modification of the technique of \cite{R1}, \cite{R2},
to control the term $I(\Phi)$. Choose $\3\in (0,2/(n-2))$. By the Jensen's 
inequality, Sobolev inequality, and (2.7),
$$\align
4\delta\tau\int_M|\nabla\Phi|^2\,dV_{\4{g}}
-\int_M\Phi^2\log\Phi^2\,dV_{\4{g}}
=&4\delta\tau\int_M|\nabla\Phi|^2\,dV_{\4{g}}
-\frac{1}{\3}\int_M\Phi^2\log\Phi^{2\3}\,dV_{\4{g}}\\
\ge&4\delta\tau\int_M|\nabla\Phi|^2\,dV_{\4{g}}
-\frac{2+2\3}{\3}\log\|\Phi\|_{L^{2+2\3}(M,\4{g})}\\
\ge&4\delta\tau\int_M|\nabla\Phi|^2\,dV_{\4{g}}
-\frac{2+2\3}{\3}\log(C\|\Phi\|_{H^1(M,\4{g})})\tag 2.9
\endalign
$$
where
$$
\|\Phi\|_{L^q(M,\4{g})}=\biggl (\int_M\Phi^q\,dV_{\4{g}}
\biggr )^{1/q}
$$
for any $q\ge 1$ and
$$
\|\Phi\|_{H^1(M,\4{g})}=\|\Phi\|_{L^2(M,\4{g})}+\|\nabla\Phi\|_{L^2(M,\4{g})}.
$$
Let
$$
\4{\Phi}=\frac{\Phi}{\|\Phi\|_{L^2(M,\4{g})}}.
$$
Then 
$$
\|\4{\Phi}\|_{H^1(M,\4{g})}\ge\|\4{\Phi}\|_{L^2(M,\4{g})}=1.\tag 2.10
$$ 
By (2.7), (2.9), and (2.10), $\forall\tau\ge\tau_0$,
$$\align
I(\Phi)\ge&4\delta\tau\|\4{\Phi}\|_{H^1(M,\4{g})}
-\frac{2+2\3}{\3}\frac{\log(C\|\4{\Phi}\|_{H^1(M,\4{g})}\|\Phi\|_{L^2(M,\4{g})})}
{\|\Phi\|_{L^2(M,\4{g})}}-4\delta\tau\\
\ge&4\delta\tau\|\4{\Phi}\|_{H^1(M,\4{g})}
-\frac{2+2\3}{\3 (4\pi\tau)^{n/2}}\log(C\|\4{\Phi}\|_{H^1(M,\4{g})})
-\frac{(1+\3)n}{\3 (4\pi\tau)^{n/2}}\log(4\pi\tau)
-4\delta\tau\\
\ge&4\delta\tau_0\|\4{\Phi}\|_{H^1(M,\4{g})}
-\frac{2+2\3}{\3 (4\pi\tau_0)^{n/2}}\log(C\|\4{\Phi}\|_{H^1(M,\4{g})})
-\frac{(1+\3)n}{\3 (4\pi\tau)^{n/2}}\log(4\pi\tau)
-4\delta\tau\\
\ge&C_{\tau_0}'-\frac{C_1}{\tau^{n/2}}\log(4\pi\tau)-4\delta\tau
\tag 2.11
\endalign
$$
where
$$
C_{\tau_0}'=\min_{y\ge 1}\biggl (4\delta\tau_0y
-\frac{2+2\3}{\3 (4\pi\tau_0)^{n/2}}\log(Cy)\biggr )>-\infty
$$
and
$$
C_1=\frac{(1+\3)n}{\3 (4\pi)^{n/2}}.
$$
By (2.8) and (2.11) we get (2.5) with $C_{\tau_0}=C_{\tau_0}'-n$. By taking
infimum over all function $f\in\Cal{A}(\4{g},\tau)$ in (2.5) we get 
$\mu(\4{g},\tau)>-\infty$ for any $\tau>0$ and the lemma follows.
\enddemo

\proclaim{\bf Corollary 2.2}
Let $\4{g}$ be a Riemannian metric on $M$. Suppose
the first eigenvalue of $R(\4{g})-4\Delta_{\4{g}}$ is positive. Then
$$
\lim_{\tau\to\infty}\mu(\4{g},\tau)=\infty.
$$
\endproclaim
\demo{Proof}
This corollary is stated without proof in \cite{P1}. We will give a
short proof of it here. We fix $\tau_0>0$ and choose $\delta\in (0,1)$
sufficiently small such that
$$
\biggl [(1-\delta)\lambda_1
-\frac{\delta}{(4\pi\tau_0)^{\frac{n}{2}}}
\|R(\4{g})\|_{L^{\infty}(M)}-4\delta\biggr ]>0.
$$
By Lemma 2.1 there exist constants $C_1>0$ and $C_{\tau_0}>0$ such that
(2.5) holds. Taking infimum over $f\in\Cal{A}(\4{g},\tau)$ in (2.5), 
we get
$$
\mu (\4{g},\tau)
\ge\biggl [(1-\delta)\lambda_1
-\frac{\delta}{(4\pi\tau_0)^{\frac{n}{2}}}
\|R(\4{g})\|_{L^{\infty}(M)}-4\delta\biggr ]\tau
-C_{\tau_0}-\frac{C_1}{\tau^{n/2}}
\log(4\pi\tau)\quad\forall\tau\ge\tau_0.
$$
Letting $\tau\to\infty$ the corollary follows.
 \enddemo

\proclaim{\bf Lemma 2.3}
Suppose $0<t_1<T$ and $\2{f}\in C^{\infty}(M\times (0,t_1))$. 
Let $\2{g}(t)=(\2{g}_{ij}(\cdot ,t))$ be an evolving metric on $M$ which 
satisfies
$$
\frac{\1}{\1 t}{\2{g}(t)}=-2(R_{ij}(\2{g}(t))+\nabla_i^{\2{g}(t)}
\nabla_j^{\2{g}(t)}\2{f})\quad\text{ in }M\times (0,t_1)\tag 2.12
$$
where $\nabla_i^{\2{g}(t)}$ is the covariant derivative with respect to
the metric $\2{g}(t)$. Suppose 
$$
\frac{\1\2{f}}{\1 t}=-\Delta_{\2{g}}\2{f}-R(\2{g})+\frac{n}{2\tau}
\quad\text{ in }M\times (0,t_1)\tag 2.13
$$
where
$$
\tau=\tau (t)=t_0'-t\tag 2.14
$$
for some constant $t_0'>t_1$. Then $\forall t\in (0,t_1)$,
$$
\frac{d}{d t}W(\2{g}(t),\2{f}(\cdot,t),\tau)
=\int_M2\tau\biggl |R_{ij}(\2{g}(t))+\nabla_i^{\2{g}(t)}\nabla_j^{\2{g}(t)}
\2{f}-\frac{1}{2\tau}\2{g}_{ij}
\biggr |^2(4\pi\tau)^{-n/2}e^{-\2{f}}\,dV_{\2{g}(t)}.
\tag 2.15
$$
\endproclaim
\demo{Proof}
This result is stated without proof in \cite{P1}. For the sake of completeness
we will give a simple proof of it here.
Let the metric $\4{g}(t)=(\4{g}_{ij}(t))$ be given by
$$
\4{g}_{ij}(t)=\frac{\2{g}_{ij}(t)}{4\pi\tau}.
$$
Then 
$$\align
W(\2{g}(t),\2{f}(\cdot,t),\tau)=&\frac{1}{4\pi}\Cal{F}(\4{g}(t),
\2{f}(\cdot,t))+\int_M(\2{f}(p,t)-n)e^{-\2{f}(p,t)}dV_{\4{g}(t)}(p)\\
\Rightarrow\quad\frac{d}{d t}W(\2{g}(t),\2{f}(\cdot,t),\tau)
=&\frac{1}{4\pi}\frac{d}{d t}\Cal{F}(\4{g}(t),\2{f}(\cdot,t))
+\int_M\2{f}_t(p,t)e^{-\2{f}(p,t)}dV_{\4{g}(t)}(p)\\
&\qquad +\int_M(\2{f}(p,t)-n)\frac{\1}{\1 t}\biggl (e^{-\2{f}(p,t)}dV_{\4{g}(t)}
(p)\biggr )\tag 2.16\endalign
$$
Now by (2.12) and (2.13),
$$\align
\frac{d}{d t}\biggl (e^{-\2{f}}dV_{\4{g}(t)}\biggr )
=&\biggl(\frac{1}{2}\4{g}^{ij}(\4{g}_{ij})_t-f_t\biggr)e^{-\2{f}}dV_{\4{g}(t)}\\
=&\biggl\{\frac{1}{2}(4\pi\tau)\2{g}^{ij}\biggl (\frac{(\2{g}_{ij})_t}
{4\pi\tau}
+\frac{\2{g}_{ij}}{4\pi\tau^2}\biggr )-\2{f}_t\biggr\}e^{-f}dV_{\4{g}(t)}\\
=&\biggl\{\2{g}^{ij}\biggl (-(R_{ij}(\2{g})
+\nabla_i^{\2{g}}\nabla_j^{\2{g}}\2{f})
+\frac{\2{g}_{ij}}{2\tau}\biggr )-\2{f}_t\biggr\}e^{-\2{f}}dV_{\4{g}(t)}\\
=&0.\tag 2.17
\endalign
$$
Hence by (2.16) and (2.17),
$$
\frac{d}{d t}W(\2{g}(t),\2{f}(\cdot,t),\tau)
=\frac{1}{4\pi}\frac{d}{d t}\Cal{F}(\4{g}(t),\2{f}(\cdot,t))+\int_M\2{f}_t(p,t)e^{-\2{f}(p,t)}
dV_{\4{g}(t)}(p).\tag 2.18
$$
By (2.17) and section 1.1 of \cite{P1},
$$\align
\frac{1}{4\pi}\frac{d}{d t}\Cal{F}(\4{g}(t),\2{f}(\cdot,t))
=&-\frac{1}{4\pi}\int_M<(\4{g}_{ij})_t,R_{ij}(\4{g})
+\nabla_i^{\2{g}}\nabla_j^{\2{g}}f>_{\4{g}}e^{-\2{f}}\,dV_{\4{g}}\\
=&-\frac{1}{4\pi}\int_M\4{g}^{ii'}\4{g}^{jj'}(\4{g}_{i'j'})_t(R_{ij}(\4{g})
+\nabla_i^{\4{g}}\nabla_j^{\4{g}}\2{f})e^{-\2{f}}\,dV_{\4{g}}\\
=&-(4\pi\tau)^{-\frac{n}{2}}\int_M\2{g}^{ii'}\2{g}^{jj'}(\tau (\2{g}_{i'j'})_t
+\2{g}_{i'j'})
(R_{ij}(\2{g})+\nabla_i^{\2{g}}\nabla_j^{\2{g}}\2{f})e^{-\2{f}}\,dV_{\2{g}}\\
=&(4\pi\tau)^{-\frac{n}{2}}\int_M[2\tau |(R_{ij}(\2{g})+\nabla_i^{\2{g}}
\nabla_j^{\2{g}}\2{f})|^2
-(R(\2{g})+\Delta_{\2{g}}\2{f})]e^{-\2{f}}\,dV_{\2{g}}\tag 2.19
\endalign
$$
By (2.13), (2.18), and (2.19), we get (2.15) and the lemma follows.
\enddemo

\proclaim{\bf Lemma 2.4}
Let $H_0\in C^{\infty}(M)$ be such that $\min_MH_0>0$. Then for any $0<t_1<T$
there exists a unique solution $H\in C^{\infty}(M\times [0,t_1])$ of the 
problem
$$\left\{\aligned
&H_t=-\Delta_{g(t)}H+R(g(t))H\quad\text{ in }M\times (0,t_1]\\
&H(x,t_1)=H_0(x)\qquad\qquad\quad\text{ in }M
\endaligned\right.\tag 2.20
$$
satisfying the condition
$$
H(x,t)\ge e^{-C_2(t_1-t)}\min_MH_0>0\quad\text{ in }M\times [0,t_1]\tag 2.21
$$
where $C_2=\|R\|_{L^{\infty}(M\times [0,t_1))}$.
\endproclaim
\demo{Proof}
By Theorem 6 of \cite{H1} there exists a unique smooth solution 
$H\in C^{\infty}(M\times [0,t_1])$ of (2.20). By continuity there exists 
$\delta_1\in (0,t_1)$ such that $H(x,t)>0$ on $M\times (t_1-\delta_1,t_1]$. 
Let
$$
t_2=\inf\{t'>0:H(x,t)>0\quad\forall x\in M,t'<t\le t_1\}.
$$
Then $0\le t'\le t_1-\delta_1$. Suppose $t'>0$. Let $s=t_1-t$. Then
$$\align
&H_s=\Delta H-R(g(t))H\ge\Delta H-C_2H\quad\forall (x,s)\in M\times 
(0,t_1-t')\\
\Rightarrow\quad&(e^{C_2s}H)_s\ge\Delta (e^{C_2s}H)
\qquad\qquad\qquad\quad\forall (x,s)\in M\times (0,t_1-t')\endalign
$$
where $C_2=\|R\|_{L^{\infty}(M\times (0,t_1))}$.
By the maximum principle for parabolic equations,
$$\align
&e^{C_2s}H\ge\min_MH_0\qquad\qquad\qquad\quad\,\,\,\forall (x,s)\in M
\times [0,t_1-t']\\
\Rightarrow\quad&H(x,t)\ge e^{-C_2(t_1-t)}\min_MH_0>0\quad\forall (x,t)\in 
M\times [t',t_1].\tag 2.22\endalign
$$
Hence by continuity there exists a constant $\delta_2\in (0,t')$ such that
$H(x,t)>0$ on $M\times (t'-\delta_2,t']$. This contradicts the maximality of
$t'$. Hence $t'=0$. Putting $t'=0$ in (2.22) we get (2.21) and the lemma 
follows.
\enddemo

\proclaim{\bf Lemma 2.5}
Let $f_0\in C^{\infty}(M)$. Then for any $0<t_1<T$, $t_0'>t_1$, there 
exists a solution $f\in C^{\infty}(M\times (0,t_1))$ of the problem
$$\left\{\aligned
&f_t=-\Delta_{g(t)}f+|\nabla f|^2-R(g(t))+\frac{n}{2\tau (t)}
\quad\text{ in }M\times (0,t_1)\\
&f(x,t_1)=f_0(x)\qquad\qquad\qquad\qquad\qquad\quad\text{ in }M
\endaligned\right.\tag 2.23
$$
where $\tau (t)=t_0'-t$.
\endproclaim
\demo{Proof}
We will use a transform of \cite{P1} to prove the lemma. Let
$$
H_0(x)=(4\pi (t_0'-t_1))^{-n/2}e^{-f_0(x)}.
$$
Then $H_0>0$ on $M$. By Lemma 2.4 there exists a unique positive solution 
$H\in C^{\infty}(M\times (0,t_1))$ of (2.20). Let
$$
f(x,t)=-\log [(4\pi\tau (t))^{n/2}H].
$$
Then by (2.20) $f$ satisfies (2.23).
\enddemo

\proclaim{\bf Theorem 2.6}
For any $t_0'>0$, $\mu (g(t),t_0'-t)$ is a monotone increasing function of 
$t\in (0,\min (t_0',T))$. If $(M,g)$ is not a Ricci soliton, then 
$\mu (g(t),t_0'-t)$ is a strictly monotone increasing function of 
$t\in (0,\min (t_0',T))$.
\endproclaim
\demo{Proof}
Let $t_1\in (0,\min (t_0',T))$. By (2.6) and an argument similar to the
proof in \cite{R1}, \cite{R2}, there exists a function $f_0\in C^{\infty}(M)$ 
satisfying
$$
(4\pi\tau(t_1))^{-n/2}\int_Me^{-f_0}dV_{g(t_1)}=1\tag 2.24
$$
and
$$
\mu (g(t_1),t_0'-t_1)=W(g(t_1),f_0,t_0'-t_1)
$$
where $\tau (t)$ is given by (2.14). Let
$f$ be the solution of (2.23) given by Lemma 2.5. Choose $t_0\in (0,t_1)$.
Let $\phi_{t_0,t}$ be as in Lemma 1.2, $\2{g}$ be given by (1.8),
and $\2{f}(p,t)=f(\phi_{t_0,t}(p),t)$. Then $\2{g}$ satisfies (1.9)
with $V(t)=-\nabla f(\cdot,t)$. 
By (2.24),
$$
(4\pi\tau(t_1))^{-n/2}\int_Me^{-\2{f}(p,t_1)}dV_{\2{g}(t_1)}=1.\tag 2.25
$$
By (1.9) and (0.1),
$$\align
\frac{\1}{\1 t}\2{g}(t)=&\phi_{t_0,t}^{\ast}\biggl (-2R_{ij}(g(t))
-2\nabla_i^{g(t)}\nabla_j^{g(t)}f\biggr )\quad\forall 0<t<t_1\\
=&-2(R_{ij}(\2{g}(t))+\nabla_i^{\2{g}(t)}\nabla_j^{\2{g}(t)}\2{f})
\quad\forall 0<t<t_1.\endalign
$$
Hence $\2{g}$ satisfies (2.12). By direct computation $\2{f}$ satisfies 
(2.13). Hence by Lemma 2.3 (2.15) holds. Thus
$$
W(\2{g}(t_1),\2{f}(\cdot ,t_1),,t_0'-t_1)\ge W(\2{g}(t),\2{f}(\cdot ,t),
,t_0'-t)+E(t,t_1)\quad\forall 0<t<t_1\tag 2.26
$$
where
$$\align
E(t,t_1)=&\int_t^{t_1}\int_M2\tau\biggl |R_{ij}(\2{g}(t))+\nabla_i^{\2{g}(t)}
\nabla_j^{\2{g}(t)}\2{f}
-\frac{1}{2\tau}\2{g}_{ij}
\biggr |^2(4\pi\tau)^{-n/2}e^{-\2{f}}\,dV_{\2{g}(t)}dt\\
=&\int_t^{t_1}\int_M2\tau\biggl |R_{ij}(g(t))+\nabla_i^{g(t)}\nabla_j^{g(t)}f
-\frac{1}{2\tau}g_{ij}
\biggr |^2(4\pi\tau)^{-n/2}e^{-f}\,dV_{g(t)}dt\\
\ge&0
\endalign
$$
with $E(t,t_1)>0$ if $g$ is not a Ricci soliton.
Since the functional $W$ is invariant under diffeomorphism, by (2.26)
$\forall 0<t<t_1$,
$$\align
&W(g(t_1),f_0,t_0'-t_1)\ge W(g(t),f(\cdot ,t),t_0'-t)+E(t,t_1)\\
\Rightarrow\quad&\mu(g(t_1),t_0'-t_1)\ge W(g(t),f(\cdot ,t),t_0'-t)+E(t,t_1).
\tag 2.27
\endalign
$$
By direct computation,
$$\align
&\frac{d}{d t}\biggl ((4\pi\tau (t))^{-n/2}\int_Me^{-f}dV_{g(t)}\biggr )=0\\
\Rightarrow\quad&(4\pi\tau(t))^{-n/2}\int_Me^{-f(p,t)}dV_{g(t)}
=(4\pi\tau(t_1))^{-n/2}\int_Me^{-f_0(p)}dV_{g(t_1)}=1\quad\forall 0<t<t_1.
\endalign
$$
Hence $f(\cdot ,t)\in\Cal{A}(g(t),\tau (t))$. Thus by (2.27),
$$
\mu (g(t_1),t_0'-t_1)\ge\mu (g(t),t_0'-t)+E(t,t_1).
$$
Since $0<t<t_1<\min (T,t_0')$ is arbitrary, the lemma follows.
\enddemo

\proclaim{\bf Theorem 2.7}
If $(M,g)$ is not a Ricci soliton, then there does not exist any shrinking 
breather for the manifold $M$ with metric $g$ evolving by the Ricci flow 
on $M\times (0,T)$.
\endproclaim
\demo{Proof}
Suppose $(M,g)$ is not a Ricci soliton and there exists a shrinking breather. 
Then there exist constants $\alpha\in (0,1)$, $0<t_1<t_2<T$, such that
$(M,\alpha g(t_1))$ is diffeomorphic to $(M,g(t_2))$. Then
$$
\mu (\alpha g(t_1),\tau)=\mu (g(t_2),\tau)\quad\forall\tau>0.\tag 2.28
$$
Since $W(\2{g},f,\tau)=W(\lambda\2{g},f,\lambda\tau)$ for any metric 
$\2{g}$ on $M$ and $\lambda>0$, for any metric 
$\2{g}$ on $M$ we have
$$
\mu (\2{g},\tau)=\mu (\lambda\2{g},\lambda\tau)\quad\forall\lambda,\tau >0.
\tag 2.29
$$
Hence by Theorem 2.6 and (2.29),
$$
\mu (\alpha g(t_1),\tau)=\mu (g(t_1),\tau/\alpha)<\mu (g(t_2),(\tau/\alpha)
-(t_2-t_1)).\tag 2.30
$$
Let $\tau=\alpha (t_2-t_1)/(1-\alpha)$. Then
$$
\frac{\tau}{\alpha}-(t_2-t_1)=\tau.\tag 2.31
$$
By (2.28), (2.30), and (2.31) we get a contradiction. Hence no shrinking 
breather exists.
\enddemo


\Refs

\ref
\key C\by B.~Chow\paper Lecture notes on Ricci flow I, II, III, 
Clay Mathematics Institute, Summer School\linebreak Program
2005 on Ricci Flow, 3-Manifolds and Geometry
June 20--July 16 at MSRI,\linebreak
http://www.claymath.org/programs/summer\_school/2005/program.php\#ricci\endref

\ref
\key CK\by \ B.~Chow and D.~Knopf\book The Ricci flow:An introduction,
Mathematical Surveys and Monographs, Volume 110, Amer. Math. Soc.
\publaddr Providence, R.I., U.S.A.\yr 2004\endref

\ref
\key H1\by R.~Hamilton\paper Three-manifolds with positive Ricci curvature
\jour J. Differential Geom.\vol 17(2)\yr 1982\pages 255--306\endref

\ref
\key H2\by R.~Hamilton\paper Four-manifolds with positive curvature
operator\jour J. Differential Geom.\vol 24(2)\yr 1986\pages 153--179\endref

\ref
\key H3\by R.S.~Hamilton\paper The Ricci flow on surfaces\jour 
Contemp. Math.\vol 71\yr 1988\pages 237--261
\endref

\ref
\key H4\by R.~Hamilton\paper The Harnack estimate for the Ricci flow
\jour J. Differential Geom.\vol 37(1)\yr 1993\pages 225--243\endref

\ref 
\key H5\by R.~Hamilton\paper The formation of singularities in the Ricci flow
\jour Surveys in differential geometry, Vol. II (Cambridge, MA, 1993),7--136,
International Press, Cambridge, MA, 1995\endref

\ref
\key Hs1\by \ \ S.Y.~Hsu\paper Global existence and uniqueness
of solutions of the Ricci flow equation\jour Differential
and Integral Equations\vol 14(3)\yr 2001\pages 305--320\endref

\ref
\key Hs2\by \ \ S.Y.~Hsu\paper Large time behaviour of solutions
of the Ricci flow equation on $R^2$\vol 197(1)\yr 2001
\pages 25--41\jour Pacific J. Math.\endref

\ref
\key Hs3\by \ \  S.Y.~Hsu\paper Asymptotic profile of
solutions of a singular diffusion equation as $t\to\infty$
\jour Nonlinear Analysis, TMA\vol 48\yr 2002\pages 781--790
\endref

\ref
\key Hs4\by \ \ S.Y.~Hsu\paper Dynamics of solutions of a
singular diffusion equation\jour Advances in Differential
Equations\vol 7(1)\yr 2002\pages 77--97\endref

\ref
\key J\by J.~Jost\book Riemannian geometry and geometric analysis, second 
edition\publ Springer-Verlag\publaddr Berlin, Heidelberg, Germany\yr 1998
\endref

\ref 
\key LY\by P.~Li and S.T.~Yau\paper On the parabolic kernel of the 
Schrodinger operator\jour Acta Math.\vol 156\yr 1986\pages 153--201\endref

\ref
\key P1\by G.~Perelman\paper The entropy formula for the Ricci flow and its 
geometric applications,\linebreak http://arXiv.org/abs/math.DG/0211159\endref 

\ref
\key P2\by G.~Perelman\paper Ricci flow with surgery on three-manifolds,
http://arXiv.org/abs/math.DG/0303109\endref

\ref
\key R1\by O.~Rothaus\paper Logarithmic Solobev inequalities and the
spectrum of Sturm-Liouville operators\jour J. Functional Analysis\vol 39
\pages 42--56\yr 1980\endref

\ref
\key R2\by O.~Rothaus\paper Logarithmic Solobev inequalities and the
spectrum of Schr\"odinger operators\jour J. Functional Analysis\vol 42
\pages 110--120\yr 1981\endref

\ref
\key S1\by W.X.~Shi\paper Deforming the metric on complete Riemannian manifolds
\jour J. Differential Geom.\vol 30\yr 1989\pages 223--301\endref

\ref
\key S2\by W.X.~Shi\paper Ricci deformation of the metric on complete 
non-compact Riemannian manifolds \jour J. Differential Geom.\vol 30\yr 1989
\pages 303--394\endref

\ref
\key W1\by \ L.F.~Wu\paper The Ricci flow on complete $R^2$
\jour Comm. in Analysis and Geometry\vol 1\yr 1993
\pages 439--472\endref

\ref
\key W2\by \ L.F.~Wu\paper A new result for the porous
medium equation\jour Bull. Amer. Math. Soc.\vol 28\yr 1993
\pages 90--94\endref

\endRefs
\enddocument